\documentclass{article}
\usepackage{amsfonts,amsmath,amssymb,cite}
\numberwithin{equation}{section}
\begin{document}
\title{Comment on the orthogonality of the Macdonald functions of
imaginary order}
\author{Rados{\l}aw Szmytkowski\footnote{Email:
radek@mif.pg.gda.pl} \hspace{0em} and
Sebastian Bielski\footnote{Email: bolo@mif.pg.gda.pl} \\*[3ex]
Atomic Physics Division,
Department of Atomic Physics and Luminescence, \\
Faculty of Applied Physics and Mathematics,
Gda{\'n}sk University of Technology, \\
Narutowicza 11/12, PL 80--233 Gda{\'n}sk, Poland}
\date{\today}
\maketitle
\begin{abstract} 
Recently, Yakubovich [Opuscula Math.\ 26 (2006) 161--172] and Passian
\emph{et al.\/} [J.\ Math.\ Anal.\ Appl.\
doi:10.1016/j.jmaa.2009.06.067] have presented alternative proofs of
an orthogonality relation obeyed by the Macdonald functions of
imaginary order. In this note, we show that the validity of that
relation may be also proved in a simpler way by applying a technique
occasionally used in mathematical physics to normalize scattering
wave functions to the Dirac delta distribution. \\*[1ex]
\end{abstract}
%
%
\section{Introduction}
\label{I}
\setcounter{equation}{0}
In recent papers, Yakubovich \cite{Yaku06} and Passian \emph{et
al.\/} \cite{Pass09} proved the following orthogonality relation for
the Macdonald functions of imaginary order:
\begin{equation}
\int\limits_{0}^{\infty}\mathrm{d}x\:
\frac{K_{\mathrm{i}\nu}(x)K_{\mathrm{i}\nu'}(x)}{x}
=\frac{\pi^{2}}{2\nu\sinh(\pi\nu)}\delta(\nu-\nu')
\qquad (\nu,\nu'>0).
\label{1.1}
\end{equation}
In the first of these works, a proof relied on advanced techniques of
the theory of distributions. An approach adopted in the second paper
was a two-step one. At first, two heuristic arguments making the
relation (\ref{1.1}) plausible were presented. One of these arguments
was based on an integral relation between the Macdonald functions of
imaginary order and Mehler's conical functions, for which a
counterpart orthogonality relation had been known for a long time.
The second argument given in support of the validity of Eq.\
(\ref{1.1}) exploited the fact that the Laplace transform of
$K_{\mathrm{i}\nu}(x)$ is a known elementary function. Subsequently,
a sophisticated proof of the relation (\ref{1.1}), different from the
one in Ref.\ \cite{Yaku06}, was presented.

It is the purpose of this note to present still another proof of the
relation (\ref{1.1}). The approach we adopt here is known in
mathematical physics, where it is occasionally used for normalization
of scattering states to the Dirac delta distribution. In some sense,
it is akin to the standard method used to prove weighted
orthogonality relations for eigenfunctions of regular
Sturm--Liouville problems (cf, e.g., Ref.\ \cite[Sec.\ 7.1]{Stak98}).
It has several advantages. First, it is \emph{elementary\/} compared
to the methods used in Refs.\ \cite{Yaku06,Pass09}. Second, it is
\emph{constructive}: one \emph{derives\/} the orthogonality relation.
Finally, it is \emph{general\/} and may be used to obtain counterpart
orthogonality relations not only for $K_{\mathrm{i}\nu}(x)$, but also
for other special functions (e.g., Ref.\ \cite{Szmy09}).
\section{A summary of relevant properties of the Macdonald functions
of imaginary order}
\label{II}
\setcounter{equation}{0}
Before we proceed to the merit, in this short section we shall
summarize these properties of the Macdonald functions of imaginary
order which will be exploited later in Sec.\ \ref{III}. The formulas
presented below have been excerpted from the collection of Magnus
\emph{et al.\/} \cite{Magn66}.

The function $K_{\mathrm{i}\nu}(x)$ is a particular solution to the
modified Bessel differential equation
\begin{equation}
x^{2}\frac{\mathrm{d}^{2}F(x)}{\mathrm{d}x^{2}}
+x\frac{\mathrm{d}F(x)}{\mathrm{d}x}+(\nu^{2}-x^{2})F(x)=0.
\label{2.1}
\end{equation}
Other particular solutions to Eq.\ (\ref{2.1}) are the modified
Bessel functions of the first kind
\begin{equation}
I_{\pm\mathrm{i}\nu}(x)
=\sum_{k=0}^{\infty}\frac{1}{k!\Gamma(k+1\pm\mathrm{i}\nu)}
\left(\frac{x}{2}\right)^{2k\pm\mathrm{i}\nu}.
\label{2.2}
\end{equation}
The relationship between the three functions is
\begin{equation}
K_{\mathrm{i}\nu}(x)=\frac{\pi}{2\mathrm{i}}
\frac{I_{-\mathrm{i}\nu}(x)-I_{\mathrm{i}\nu}(x)}{\sinh(\pi\nu)},
\label{2.3}
\end{equation}
from which it follows immediately that
\begin{equation}
K_{\mathrm{i}\nu}(x)=K_{-\mathrm{i}\nu}(x).
\label{2.4}
\end{equation}
For large positive values of $x$, the function $K_{\mathrm{i}\nu}(x)$
has the asymptotic representation
\begin{equation}
K_{\mathrm{i}\nu}(x)\stackrel{x\to\infty}{\sim}
\sqrt{\frac{\pi}{2x}}\mathrm{e}^{-x}[1+O(x^{-1})],
\label{2.5}
\end{equation}
while for $x\to+0$ from Eqs.\ (\ref{2.3}) and (\ref{2.2}) and with
the aid of the known relationship
\begin{equation}
|\Gamma(\mathrm{i}\nu)|=\sqrt{\frac{\pi}{\nu\sinh(\pi\nu)}}
\qquad (\nu\in\mathbb{R})
\label{2.6}
\end{equation}
one deduces that
\begin{eqnarray}
&& K_{\mathrm{i}\nu}(x) \stackrel{x\to+0}{\sim}
\sqrt{\frac{\pi}{\nu\sinh(\pi\nu)}}
\cos\left[-\nu\ln\frac{x}{2}+\arg\Gamma(\mathrm{i}\nu)\right]
+O\left(x^{2}\sin\left[-\nu\ln\frac{x}{2}
+\arg\Gamma(2+\mathrm{i}\nu)\right]\right)
\nonumber \\
&& \hspace*{30em} (\nu\in\mathbb{R}).
\label{2.7}
\end{eqnarray}
\section{Derivation of the orthogonality relation for the Macdonald
functions of imaginary order}
\label{III}
\setcounter{equation}{0}
To derive the orthogonality relation for the Macdonald functions of
imaginary order, we proceed as follows. If $K_{\mathrm{i}\nu}(x)$,
with $\nu\in\mathbb{R}$, is substituted for $F(x)$ into Eq.\
(\ref{2.1}), this results in the differential identity
\begin{equation}
\frac{\mathrm{d}}{\mathrm{d}x}
\left(x\frac{\mathrm{d}K_{\mathrm{i}\nu}(x)}{\mathrm{d}x}\right)
+\left(\frac{\nu^{2}}{x}-x\right)K_{\mathrm{i}\nu}(x)=0.
\label{3.1}
\end{equation}
The counterpart identity for the function $K_{\mathrm{i}\nu'}(x)$,
with $\nu'\in\mathbb{R}$, is
\begin{equation}
\frac{\mathrm{d}}{\mathrm{d}x}
\left(x\frac{\mathrm{d}K_{\mathrm{i}\nu'}(x)}{\mathrm{d}x}\right)
+\left(\frac{\nu^{\prime\,2}}{x}-x\right)K_{\mathrm{i}\nu'}(x)=0.
\label{3.2}
\end{equation}
Next, we premultiply the first of the above equations by
$K_{\mathrm{i}\nu'}(x)$, the second one by $K_{\mathrm{i}\nu}(x)$,
subtract and integrate the result over $x$ from $x=\xi>0$ to
$x=\infty$. After obvious rearrangements, this gives
\begin{equation}
(\nu^{2}-\nu^{\prime\,2})\int\limits_{\xi}^{\infty}\mathrm{d}x\:
\frac{K_{\mathrm{i}\nu}(x)K_{\mathrm{i}\nu'}(x)}{x}
=\int\limits_{\xi}^{\infty}\mathrm{d}x\:
\left[K_{\mathrm{i}\nu}(x)\frac{\mathrm{d}}{\mathrm{d}x}
\left(x\frac{\mathrm{d}K_{\mathrm{i}\nu'}(x)}{\mathrm{d}x}\right)
-K_{\mathrm{i}\nu'}(x)\frac{\mathrm{d}}{\mathrm{d}x}
\left(x\frac{\mathrm{d}K_{\mathrm{i}\nu}(x)}
{\mathrm{d}x}\right)\right].
\label{3.3}
\end{equation}
The integral on the right-hand side of Eq.\ (\ref{3.3}) is easily
evaluated by parts; one obtains
\begin{equation}
(\nu^{2}-\nu^{\prime\,2})\int\limits_{\xi}^{\infty}\mathrm{d}x\:
\frac{K_{\mathrm{i}\nu}(x)K_{\mathrm{i}\nu'}(x)}{x}
=\left[x\left(K_{\mathrm{i}\nu}(x)
\frac{\mathrm{d}K_{\mathrm{i}\nu'}(x)}{\mathrm{d}x}
-K_{\mathrm{i}\nu'}(x)\frac{\mathrm{d}K_{\mathrm{i}\nu}(x)}
{\mathrm{d}x}\right)\right]_{x=\xi}^{\infty}.
\label{3.4}
\end{equation}
By virtue of Eq.\ (\ref{2.5}), the expression in the bracket on the
right-hand side of the above relation vanishes in the upper limit.
Hence, we obtain
\begin{equation}
\int\limits_{\xi}^{\infty}\mathrm{d}x\:
\frac{K_{\mathrm{i}\nu}(x)K_{\mathrm{i}\nu'}(x)}{x}
=-\,\xi\frac{\displaystyle K_{\mathrm{i}\nu}(\xi)
\frac{\mathrm{d}K_{\mathrm{i}\nu'}(\xi)}{\mathrm{d}\xi}
-K_{\mathrm{i}\nu'}(\xi)\frac{\mathrm{d}K_{\mathrm{i}\nu}(\xi)}
{\mathrm{d}\xi}}{\nu^{2}-\nu^{\prime\,2}}
\label{3.5}
\end{equation}
and consequently
\begin{equation}
\int\limits_{0}^{\infty}\mathrm{d}x\:
\frac{K_{\mathrm{i}\nu}(x)K_{\mathrm{i}\nu'}(x)}{x}
=-\lim_{\xi\to+0}\xi\frac{\displaystyle K_{\mathrm{i}\nu}(\xi)
\frac{\mathrm{d}K_{\mathrm{i}\nu'}(\xi)}{\mathrm{d}\xi}
-K_{\mathrm{i}\nu'}(\xi)\frac{\mathrm{d}K_{\mathrm{i}\nu}(\xi)}
{\mathrm{d}\xi}}{\nu^{2}-\nu^{\prime\,2}}.
\label{3.6}
\end{equation}
Use of the asymptotic representation (\ref{2.7}) and of elementary
trigonometric identities transforms Eq.\ (\ref{3.6}) into
\begin{eqnarray}
\int\limits_{0}^{\infty}\mathrm{d}x\:
\frac{K_{\mathrm{i}\nu}(x)K_{\mathrm{i}\nu'}(x)}{x}
&=& \frac{\pi}{2\sqrt{\nu\nu'\sinh(\pi\nu)\sinh(\pi\nu')}}
\nonumber \\
&& \times\lim_{\xi\to+0}\left\{\frac{\displaystyle
\sin\left[-(\nu-\nu')\ln\frac{\xi}{2}+\arg\Gamma(\mathrm{i}\nu)
-\arg\Gamma(\mathrm{i}\nu')\right]}{\nu-\nu'}\right.
\nonumber \\
&& \left.+\,\frac{\displaystyle
\sin\left[-(\nu+\nu')\ln\frac{\xi}{2}+\arg\Gamma(\mathrm{i}\nu)
+\arg\Gamma(\mathrm{i}\nu')\right]}{\nu+\nu'}\right\}.
\label{3.7}
\end{eqnarray}
To evaluate the limit on the right-hand side of Eq.\ (\ref{3.7}), we
observe that if $f(\eta)$ is a real analytic function of
$\eta\in\mathbb{R}$, such that $f(0)=0$ (which implies that
$\lim_{\eta\to0}f(\eta)/\eta$ is finite), then in the distributional
sense it holds that
\begin{equation}
\lim_{a\to\infty}\frac{\sin[a\eta+f(\eta)]}{\pi\eta}
=\frac{1}{2\pi}\lim_{a\to\infty}
\int\limits_{-a-f(\eta)/\eta}^{a+f(\eta)/\eta}\mathrm{d}\alpha\:
\mathrm{e}^{\mathrm{i}\alpha\eta}
=\frac{1}{2\pi}\lim_{a\to\infty}
\int\limits_{-a}^{a}\mathrm{d}\alpha\:
\mathrm{e}^{\mathrm{i}\alpha\eta}
=\frac{1}{2\pi}\int\limits_{-\infty}^{\infty}\mathrm{d}\alpha\:
\mathrm{e}^{\mathrm{i}\alpha\eta}.
\label{3.8}
\end{equation}
In the expression at the extreme right of the above chain of
equalities one immediately recognizes the well-known Fourier
representation of the Dirac delta distribution $\delta(\eta)$ (cf,
e.g., \cite[Sec.\ 4.5]{Sned51}), so that, provided the function
$f(\eta)$ satisfies the above constraints, one has
\begin{equation}
\lim_{a\to\infty}\frac{\sin[a\eta+f(\eta)]}{\pi\eta}=\delta(\eta).
\label{3.9}
\end{equation}
As for $\xi\to+0$ it holds that $-\ln(\xi/2)\to\infty$, with the help
of the above relationship Eq.\ (\ref{3.7}) becomes
\begin{equation}
\int\limits_{0}^{\infty}\mathrm{d}x\:
\frac{K_{\mathrm{i}\nu}(x)K_{\mathrm{i}\nu'}(x)}{x}
=\frac{\pi^{2}}{2\sqrt{\nu\nu'\sinh(\pi\nu)\sinh(\pi\nu')}}
[\delta(\nu-\nu')+\delta(\nu+\nu')].
\label{3.10}
\end{equation}
Exploiting in Eq.\ (\ref{3.10}) the following basic property of the
delta distribution \cite[Sec.\ 4.4]{Sned51}:
\begin{equation}
g(\eta')\delta(\eta-\eta')=g(\eta)\delta(\eta-\eta'),
\label{3.11}
\end{equation}
one eventually arrives at the sought orthogonality relation
\begin{equation}
\int\limits_{0}^{\infty}\mathrm{d}x\:
\frac{K_{\mathrm{i}\nu}(x)K_{\mathrm{i}\nu'}(x)}{x}
=\frac{\pi^{2}}{2\nu\sinh(\pi\nu)}
[\delta(\nu-\nu')+\delta(\nu+\nu')].
\label{3.12}
\end{equation}
If we impose the constraint $\nu,\nu'>0$, then $\nu+\nu'>0$ and
consequently in the distributional sense we have
\begin{equation}
\delta(\nu+\nu')=0
\qquad (\nu,\nu'>0).
\label{3.13}
\end{equation}
It is then evident that under the above restriction the orthogonality
relation (\ref{3.12}) turns into the one in Eq.\ (\ref{1.1}).
%
%

%
\end{document}